\input amssym.def
\input amssym
\magnification=1200
\parindent0pt
\hsize=16 true cm
\baselineskip=13  pt plus .2pt
$ $

\def\Z{\Bbb Z}
\def\A{\Bbb A}
\def\S{\Bbb S}

\centerline {\bf On finite simple and nonsolvable groups
acting on homology 4-spheres}

 \bigskip \bigskip

\centerline {Mattia Mecchia and  Bruno Zimmermann}

\bigskip

\centerline {Universit\`a degli Studi di Trieste}
\centerline {Dipartimento di Matematica e Informatica}
\centerline {34100 Trieste, Italy}

\bigskip \bigskip

Abstract.  {\sl  The only finite nonabelian simple group acting on a
homology 3-sphere - necessarily non-freely - is the dodecahedral group
$\Bbb A_5 \cong {\rm PSL}(2,5)$ (in analogy, the only
finite perfect group acting freely on a homology 3-sphere is the binary
dodecahedral group  $\Bbb A_5^* \cong {\rm SL}(2,5)$). In the present
paper we show that the only finite simple groups acting on a homology
4-sphere, and in particular on the 4-sphere, are the alternating or
linear fractional groups groups
$\Bbb A_5 \cong {\rm PSL}(2,5)$ and
$\Bbb A_6 \cong {\rm PSL}(2,9)$. From this we deduce a short list of groups
which contains all finite nonsolvable groups
admitting an action on a homology 4-spheres. }


\smallskip


\bigskip \bigskip

{\bf 1. Introduction}

\medskip

We are interested in finite groups, and in particular in  finite
simple and nonsolvable groups admitting  smooth
orientation-preserving actions on integer homology spheres (in
the present paper, simple group will always mean  {\it
nonabelian} simple group; also, all actions will be smooth and
orientation-preserving). Let $G$ be a finite group acting on a
homology $n$-sphere. If the action is free than $G$ has periodic
cohomology, of period dividing $n+1$, and the groups of periodic
cohomology have been classified by Zassenhaus and Suzuki; the
only perfect groups among them are the linear groups
${\rm SL}(2,p)$, for prime numbers $p\ge 5$. The cohomological
period of
${\rm SL}(2,p)$ is the least common multiple of 4 and $p-1$
(see [Sj],[Sw]), so the only finite perfect group acting freely
on a homology 3-sphere is the binary dodecahedral group
$\A_5^* \cong {\rm SL}(2,5)$ (see also [Mn]).

\medskip

In the present paper we are interested in arbitray, i.e. not
necessarily free actions of finite nonsolvable groups on homology
spheres. We note that a finite nonabelian simple group
does not admit free actions on a homology sphere (any simple group has
a subgroup $\Z_2 \times \Z_2$ which does not act freely). It is shown
in [Z1] that the only finite simple group acting on a homology 3-sphere
is the dodecahedral group $\A_5 \cong {\rm PSL}(2,5)$, and the
finite nonsolvable groups acting on a homology 3-sphere are determined
in [MeZ] and [Z2] and are closely related to either $\A_5$ or $\A_5^*$.
In the present paper, we determine the finite simple and nonsolvable
groups  which admit an orientation-preserving action on a homology
4-sphere (necessarily non-free by the Lefschetz fixed point
theorem).  Our first main result is the following.

\bigskip

{\bf 1.1 Theorem}  {\sl The only finite  nonabelian simple
groups  admitting an action on a homology 4-sphere, and in particular
on the 4-sphere, are the alternating or linear fractional groups $\A_5
\cong {\rm PSL}(2,5)$ and  $\A_6 \cong {\rm PSL}(2,9)$.}

\bigskip

We note that the alternating group $\A_6$ acts on the 5-simplex
(permuting its vertices), and hence on its boundary
homeomorphic to the 4-sphere.

\bigskip

Using Theorem 1.1 we will prove then the following.

\bigskip

{\bf 1.2 Theorem} {\sl  Let $G$ be a finite non-solvable group
acting orientation-preservingly on a homology 4-sphere. Then one
of the following cases occurs:

\medskip

a) $G$ contains a normal subgroup  isomorphic to $(\Bbb{Z}_2)^4$,
with factor group isomorphic  to  $\Bbb{A}_5$ or  the
symmetric group $\Bbb{S}_5$;

\medskip

b) $G$ is isomorphic to $\A_6$ or $\S_6$;

\medskip

c) $G$ contains, of index at most two,  a  subgroup isomorphic
to  one of the following groups:

- $\Bbb{A}_5\times C$  where  $C$ is dihedral or cyclic;

- the central product $\Bbb{A}_5^*\times_{\Bbb{Z}_2}\Bbb{A}_5^*$;

- a central product $\Bbb{A}_5^* \times_{\Bbb{Z}_2} C$ where  $C$ is a
solvable group that  admits a  free  and orientation-preserving
action on the 3-sphere.}

\bigskip

We remark that the groups listed in part c) of the Theorem are close
to the class of nonsolvable groups admitting an action already in
dimension three, that is on a homology 3-sphere, see [MeZ] and [Z2]
for lists of such groups. In contrast, the groups in a) and b)
do not admit an action on a homology
3-sphere.  The natural candidate for a group of type a) is obtained
from the semidirect product $(\Z_2)^5 \ltimes \S_5$ where
$\S_5$ acts on the normal subgroup $(\Z_2)^5$ by permuting the
components (Weyl group or wreath product $\Z_2 \wr \S_5$). The group
$(\Z_2)^5 \ltimes \S_5$ acts orthogonally on euclidean 5-space by
inversion and permutation of coordinates, and the subgroup of index
two of orientation-preserving elements is a semidirect product
$(\Z_2)^4  \ltimes \S_5$ which acts orthogonally on the 4-sphere.
Concerning case b), the symmetric group $\S_6$ acts on the 5-simpex by
permuting its vertices, and on its boundary which is the 4-sphere; by
composing the orientation-reversing elements with
$-id_{\S^4}$ one obtains an orientation-preserving action of $\S_6$
on $S^4$.

\bigskip

We close the introduction with some remarks on the situation in
higher dimensions.  There is no finite perfect group
acting freely on a homology 5-sphere (because none of the groups
${\rm SL}(2,p)$ has period six: after
${\rm SL}(2,5)$, of period four, the next cases are
${\rm SL}(2,7)$ and ${\rm SL}(2,13)$, of period 12). Considering
the case of simple groups, besides $\A_5$ and $\A_6$ there are
at least three other finite simple groups acting on a homology
5-sphere which are $\A_7$, its subgroup ${\rm PSL}(2,7)$, and the
unitary group $U_4(2)$ (in the notation of [A]; this is a
subgroup of index two in the Weyl or Coxeter group of type $E_6$
which has an 6-dimensional integral linear representation); all
these groups act orthogonally on the 5-sphere. In
dimension six, there are orthogonal actions on the 6-sphere of
the simple groups $\A_8$, ${\rm PSL}(2,8)$, ${\rm PSL}(2,13)$,
$U_3(3)$ and of the symplectic group $S_6(2)$ (a subgroup of
index two in the Weyl group of type $E_7$). We do not know if
there are still other simple groups acting in these dimensions
(but we suppose that the answer is no). This naturally leads to the
following

\bigskip

{\bf Problems} i)  What is the minimal dimension of a homology
sphere on which ${\rm PSL}(2,q)$ acts? (or a metacyclic group
$H(p:q)$, see Lemma 2.2; there should be some analogy here with the
case of free actions of the groups ${\rm SL}(2,p)$ where a lower bound
comes from the cohomological dimension)

ii) Show that in every dimension there are only finitely many
groups ${\rm PSL}(2,q)$, and more generally only finitely many finite
simple groups acting on an integer homology sphere.

iii) Show that every group ${\rm PSL}(2,q)$, with $q$ an odd prime
power,  acts on a $\Z_2$-homology 3-sphere (i.e., with coefficients in
the integers mod two, see [MeZ],[Z3])

\bigskip  \bigskip

{\bf 2. Preliminary results}

\medskip

The first step of the proof of Theorem 1.1 is the following

\bigskip

{\bf 2.1 Proposition} {\sl  Let $G$ be a finite group acting
orientation-preservingly on a homology 4-sphere.

a) If $G$ is a linear fractional group ${\rm PSL}(2,q)$ then $q
\le 5$ or $q=9$.

b) If $G$ is a linear group ${\rm SL}(2,q)$ then  $q \le 5$.}

\bigskip

We start with the following two Lemmas.

\bigskip

{\bf 2.2 Lemma} {\sl  For a prime $p$ and an integer $q\ge 2$,
let $H=H(p:q)$ be a metacyclic group, with normal subgroup $\Z_p$
and factor group $\Z_q$, acting orientation-preservingly on a
homology $m$-sphere.

a) If $m=3$ then any element of $\Z_q$ acts by  $\pm$identity on
$\Z_p$.

b) If $m=4$ then the square of any element of $\Z_q$ acts by
$\pm$identity on $\Z_p$.}

\bigskip

{\it Proof.}  a) Follows from [Z1,proof of Proposition 1].

b) By [Bo,chapter IV.4], the fixed point set of the normal
subgroup
$\Z_p$ of $H$ is a homology sphere of even codimension, so the
fixed point set is a 0-sphere or a 2-sphere.

\medskip

Suppose first that the fixed point set of $\Z_p$ is a 2-sphere
$S^2$. Then,  locally in 4-dimensional space, $\Z_p$ acts as
rotations around
$S^2$; also,
$S^2$ is invariant under the action of $H$. Any element of $H$
conjugates a rotation in $\Z_p$ of minimal angle around $S^2$ to
a rotation of minimal angle, and hence induces $\pm$identity on
$\Z_p$ (see [Br,chapter VI.2] for the existence of equivariant
tubular neighbourhoods).

\medskip

Suppose then that the fixed points set of $\Z_p$ is a 0-sphere
$S^0$. Again $S^0$ is invariant under the action of $H$, and a
subgroup $H_0$ of index one or two in $H$ fixes both points of
$S^0$. The boundary of an
$H_0$-invariant neighbourhood of one of these two fixed points is
homeomorphic to a 3-sphere, so $H_0$ acts on a homology
3-sphere. Now b) follows from a).

\bigskip

{\bf 2.3 Lemma}  {\sl For a prime $p$, let $A$ be an elementary
abelian
$p$-group acting orientation-preservingly on a homology m-sphere.

a) If $m=3$ then $A$ has rank at most two if $p>2$, and rank at
most three if $p=2$.

b) If $m=4$ then $A$ has rank at most two if $p>2$, and rank at
most four if $p=2$.}

\bigskip

{\it Proof.}  a) See [MeZ,Proposition 4].

b)  Consider a subgroup $\Z_p$ of
$A$. The fixed points set of $\Z_p$ is a 0-sphere $S^0$ or a
2-sphere
$S^2$, invariant under the action of $A$.

\medskip

Suppose first that the fixed point set of $\Z_p$ is $S^2$. Any
finite orientation-preserving group fixing $S^2$ pointwise acts
as rotations around
$S^2$ and hence is cyclic; this implies that the factor group
$A/\Z_p$ acts faithfully on $S^2$. The elementary abelian
$p$-groups acting on a 2-sphere are cyclic if $p>2$, and of rank
at most three if $p=2$. Thus $A$ is as stated in the Lemma.

\medskip

Now suppose that the fixed point set of $\Z_p$ is a 0-sphere
$S^0$. Then a subgroup $A_0$ of index one or two fixes both
points of $S^0$, and $A_0$ maps a 3-sphere to itself which is
the boundary of a regular invariant neighbourhood of one of  the
two fixed points. By  a), the only elementary abelian group of
rank $\ge 3$ acting orientation-preservingly on a homology
3-sphere is the group $(\Z_2)^3$ which implies the statement of
the Lemma.

\bigskip

{\it Proof of Proposition 2.1}

\medskip

a) Suppose that  $G = {\rm PSL}(2,q)$ acts on a homology 3-sphere
$M$, for a prime power $q=p^n$. We will show that $q \le 5$ or
$q = 9$.

\medskip

Let $A \cong (\Z_p)^n$ be the subgroup of all elements of $G$
represented by upper triangular matrices of the form
$$  \pmatrix {
                 1  &  y  \cr
                 0  &  1  \cr},$$ where $y \in GF(q)$ (so $A$ is
isomorphic to the additive group of the the finite Galois field
$GF(q)$). For a generator $x$ of the multiplicative group of the
field
$GF(q)$, let $B$ be the cyclic subgroup of $G$ generated by the
matrix
$$  \pmatrix {
                   x  &    0     \cr
                   0  &  x^{-1}  \cr},$$

of order $(q-1)/2$ if $p>2$, and of order $q-1$ if $p=2$. The
cyclic group $B$ normalizes $A$, so the subgroup $H$ generated
by $A$ and $B$ is a semidirect product of these two groups
(metacyclic if $A$ is cyclic). It is easy to see that no
nontrivial element of $B$ operates trivially on $A$ by
conjugation, or in other words that the action of $B$ on $A$ is
faithful.

\medskip

Suppose first that $p$ is an odd prime. For the subgroup
${\rm PSL}(2,p)$ of ${\rm PSL}(2,p^n)$, the subgroup $A$ is
cyclic of order $p$ and $B$ cyclic of order $(p-1)/2$, so Lemma
2.2 implies
$p=3$ or $p=5$.  For $G = {\rm PSL}(2,p^n)$, we have $A \cong
(\Z_p)^n$, and Lemma 2.3 implies $n \le 2$. This leaves us with
the two groups
${\rm PSL}(2,9)$ (which admits an action on the 4-sphere) and
${\rm PSL}(2,25)$.

\medskip

The group ${\rm PSL}(2,25)$ can be excluded by applying the Borel
formula [Bo,Theorem XIII.2.3] to the subgroup $A \cong (\Z_5)^2$. The
Borel formula states that

$$4-r = \Sigma_H (n(H)-r)$$

where the sum is taken over the six different subgroups
$H \cong \Z_5$ of $A$, $n(H)$ denotes the
dimension of the fixed point set of a subgroup $H$ and $r$ the
dimension of the fixed point set of $A \cong (\Z_5)^2$ (equal to
-1 if the fixed point set is empty). We note that all elements
of order five in ${\rm PSL}(2,25)$ are conjugate, so the formula
becomes $4 + 5r = 6 n(H)$ which admits no integer solution $n(H)$
for $-1 \le r \le 3$ and hence excludes ${\rm PSL}(2,25)$.

\medskip

Now suppose that $q=2^n$ is a power of two.
Then ${\rm PSL}(2,q)$ has a subgroup $A \cong (\Z_2)^n$ and Lemma
2.3 implies $n \le 4$. The groups
${\rm PSL}(2,2) \cong S_3$ and ${\rm PSL}(2,4) \cong \A_5$ act on
the 4-sphere which leaves us with the groups ${\rm PSL}(2,8)$
and ${\rm PSL}(2,16)$.

\medskip

We apply again the Borel formula to the subgroup
$A \cong (\Z_2)^3$ of ${\rm PSL}(2,8)$ where the sum is now taken
over the seven subgroups $H \cong (\Z_2)^2$ of index two in $A$.
The subgroup $A \cong (\Z_2)^3$ of ${\rm PSL}(2,8)$  is
normalized by a cyclic subgroup $B$ of order seven which acts
transitively on the seven involutions of $A$; then $B$ acts
transitively also on the seven subgroups $H \cong (\Z_2)^2$ of
index two in $A$, and the Borel formula becomes $4+6r = 7n(H)$;
again this has no solution for $-1 \le r \le 3$, so ${\rm
PSL}(2,8)$ does not act.

\medskip

Finally, we consider ${\rm PSL}(2,16)$ and its subgroup $A \cong
(\Z_2)^4$.  The group $B$ is cyclic of order 15 now and acts
transitively on the 15 involutions in  $A$ which hence are all
conjugate. Then also all 15  subgroups  $H \cong (\Z_2)^3$ of
index two in $A$ are conjugate, and the Borel formula reads
$4+14r = 15n(H)$ which again has no solution.  This excludes also
${\rm PSL}(2,16)$ and finishes the proof of part a) of the
Proposition.

\medskip

b) Suppose that  $G = {\rm SL}(2,q)$ acts on a homology 4-sphere
$M$, for a prime power $q=p^n$. By [Bo,chapter IV.4], the central
involution $z$ of
${\rm SL}(2,q)$ has fixed point set $S^0$ or $S^2$. Suppose
first that it has fixed point set
$S^2$. Then $S^2$ is invariant under the action of $G$, and $G$
induces an  action of ${\rm SL}(2,q)/z \cong {\rm PSL}(2,q)$ on
$S^2$ which implies $q \le 5$. Now suppose that $z$ has fixed
point set $S^0$. Then
$G$ fixes each of the two points in $S^0$ and acts on a regular
neighbourhood of each of them which is a 3-sphere. By [Z2] or
[MeZ], the only perfect group of type ${\rm SL}(2,q)$ acting on a
homology 3-sphere is the binary dodecahedral group
$\A_5^* \cong {\rm SL}(2,5)$, so $q \le 5$.

\medskip

This finishes the proof of Proposition 2.1. The proof of part b)
generalizes to give the following Lemma which will be used in
the proof of Theorem 1.2 (the 3-dimensional case follows from
[Z2] or [MeZ]).

\bigskip

{\bf 2.4 Lemma}  {\sl Let $G$ be a central extension, with
nontrivial center, of a nonabelian simple group $\bar G$. If
$G$ acts on a homology 3- or  4-sphere then $\bar G$ is
isomorphic to the dodecahedral group $\A_5 \cong {\rm
PSL}(2,5)$. }

\bigskip

{\bf 3. Proof of Theorem 1.1}

\medskip

We will show that a finite simple group acting on a homology
4-sphere has sectional 2-rank at most four (i.e., every 2-subgroup is
generated by at most four elements), and then apply the
Gorenstein-Harada classification of such groups. Then, by suitable
subgroup considerations using Proposition 2.1 and Lemmas 2.2 - 2.4, all
groups of the Gorenstein-Harada list can be excluded except
$\A_5$ and $\
A_6$ (we note that this works also
for the full list of the finite simple groups).

\bigskip

{\bf 3.1 Proposition}

\smallskip

{\sl Let $S$ be a finite 2-group acting orientation-preservingly on
a $\Z_2$-homology 4-sphere $M$. Then $S$ is generated by at most four
elements (i.e., has rank at most four).}

\bigskip

{\it Proof.}  Let $h$ be a central involution in $S$. By Smith fixed
point theory ([Br]), the fixed point set of Fix($h$) of $h$ is a
2-sphere $S^2$ or a 0-sphere $S^0$ (since $h$ is orientation
preserving, the codimension of the fixed point set is even); since
$h$ is central, Fix($h$) is invariant under the action of $S$.

\medskip

Suppose first that Fix($h$) is a 2-sphere $S^2$. Taking an
equivariant tubular neighbourhood, the subgroup
$F$ of $S$ which fixes $S^2$ pointwise acts orientation-preservingly
on an orthogonal 1-sphere $S^1$ and hence is cyclic. The factor
group
$S/F$ acts effectively on $S^2$; such an action is conjugate to an
orthogonal one, in particular $S/F$ is isomorphic to a finite
subgroup of the orthogonal group O(3). The finite subgroups of O(3)
are well-known (cyclic, dihedral, tetrahderal, octahedral,
dodecahedral groups and 2-fold extensions of such groups); these
groups have rank at most three, and hence $S$ has rank  at most
four.

\medskip

Now suppose that Fix($h$) is a 0-sphere $S^0$ (i.e., two points). By
[DH] the finite 2-group $S$ admits also an orthogonal action on the
4-sphere $S^4$ such that the dimension of the fixed point set of
any subgroup of $S$ coincides for the actions of $S$ on $M$ and
$S^4$  (the two actions have the same dimension function for the
fixed point set of subgroups of $S$).  We consider such an
orthogonal action of
$S$ on $S^4$. The involution $h$ fixes a 0-sphere $S^0$ which
consists of two antipodal points of $S^4$, and we denote $S^3$ the
corresponding equatorial 3-sphere between these two antipodal
points. Note that with $S^0$ also $S^3$ is invariant under the
action of $S$, and hence $S$ is isomorphic to a subgroup of the
orthogonal group O(4) of $S^3$. A list of the finite subgroups  of
O(4) can be found in [DV].  One way to proceed now is to identify
the finite 2-groups among these groups and prove one by one that
they are generated by at most four elements. Since the lists in
[DV] are rather technical and a proof of completeness is basically
not given, we prefer to present a direct argument avoiding such a
list.

\medskip

We can assume that $S$ is a finite subgroup of the orthogonal group
O(4); we will prove that $S$ has rank at most four.

\medskip

If every abelian normal subgroup of $S$ is cyclic then, by
[Su2,4.4.3], $S$ is a cyclic,
dihedral, quaternion or quasi-dihedral 2-group; since each of these
groups has rank at most two, this finishes the proof in this case.

\medskip

We can assume then that $S$ has a non-cyclic abelian normal
subgroup. By [Su2,4.4.5], $S$ has also a normal subgroup
$U = \Z_2 \times \Z_2$. The group $S$ acts by conjugation on the
three involutions of $U$; since $S$ is a finite 2-group, each
involution in $U$ is fixed by a subgroup of index one or two in $S$.
If some involution in $U$ has fixed point set $S^2$ then $S^2$ is
invariant under $S$ or a subgroup of index two; since every finite
2-subgroup of O(3) has rank at most three, $S$ has rank at most
four. Similarly, if some involution in $U$ has fixed point set
$S^0$ (two antipodal points of $S^3$) then again $S$ or a subgroup
of index two leaves invariant $S^0$ and the corresponding
equatorial 2-sphere
$S^2$, so again $S$ has rank at most four.

\medskip

Since $U = \Z_2 \times \Z_2$ does not act freely on $S^3$, we can
assume then that some involution  $u$ in $U$ has fixed point set
$S^1$. If $u$ is central in $S$ then $F$ is invariant under $S$ and
it is again easy to see that $S$ has rank at most four. Now the
remaining case is the following: two involutions $u_1$ and
$u_2$ in $U$ have fixed point set $F_1 \cong S^1$ and $F_2 \cong
S^1$ and are exchanged by some element of $S$, their product is
central in
$S$ and acts freely. Note that $F_1$ and $F_2$ are great circles in
$S^3$; decomposing $\Bbb R^4$ as $\Bbb R^2 \times
\Bbb R^2$, these are obtained by intersecting $S^3$ with the two
orthogonal planes $\Bbb R^2$ of such a decomposition. Denoting by
$D_n$ a rotation of $\Bbb R^2$ of order $n$, one has $u_1 =
(D_1,D_2)$ and
$u_2 = (D_2,D_1)$.

\medskip

Let $S_0 = S \; \cap $ SO(4) be the orientation-preserving subgroup,
of index one or two in $S$, with
$U \subset S_0$.  An element of $S$ either exchanges $F_1$ and
$F_2$, or leaves invariant both $F_1$ and $F_2$; if an element of
$S_0$ maps
$F_1$ to itself then it  acts as a rotation around $F_1$ (fixing
$F_1$ pointwise), or as a rotation along $F_1$, or as strong
inversion of $F_1$. Let $A$ be the subgroup of elements of $S_0$
which are rotations of $F_1$. Then $A$ contains exactly the three
involutions $u_1,u_2$ and $u_1u_2$. Every element  of $A$ acts as a
rotation also on $F_2$, and a strong inversion of $F_1$ is also a
strong inversion of $F_2$.  Note that
$A$ is normal in $S$ and a direct product of two nontrivial cyclic
groups.  If there are no strong inversions or no elements in $S_0$
which exchage $F_1$ and $F_2$ then $S_0$ has rank at most three and
we are done. So we can assume that there is a strong inversion $r$
of
$F_1$, and also an element
$s$ in $S_0$ exchanging $F_1$ and $F_2$. Note that $r$ and $s$
together with $A$ generate
$S_0$. If $S=S_0$ then $S$ has rank at most four and we are done, so
we can assume that there is also an orientation-reversing element
$t$ in $S$. By eventually composing
$t$ with $s$, we can assume that $t$ leaves invariant both $F_1$ and
$F_2$.  By eventually composing $t$ with $r$, we can assume that $t
= (D,R)$ where $D$ denotes a rotation and $S$  a reflection of
$\Bbb R^2$.

\medskip

For powers $x,y$ and $z$ of two, we denote by $h = (D_1,D_x)$ a
rotation of maximal order (or minimal angle) around $F_1$ (i.e.,
with fixed point set $F_1$), and by
$k = (D_y,D_z)$ a rotation of minimal translation length  along
$F_1$ (so $k^y$ is the minimal power of $k$ which fixes
$F_1$ pointwise). Note that $h$ and
$k$ generate $A$. If $x \ge z$ then we can choose
$k = (D_y,D_1)$; since $s$ exchanges $F_1$ and $F_2$ this implies
$x=y$.  In this case $s$ and $h$ generate $k$, and $h,r,s$ and $t$
generate $S$.

\medskip

So we can assume that $x<z$. Since $A$ is normal in $S$, we have
$tkt^{-1} = (D_y,D_z^{-1}) = k(D_1,D_z^{-2}) = kh^q$, for a
primitive power $h^q$ of $h$. Then $t$ and $k$ generate
$h$, and  $k,r,s$ and $t$ generate $S$.

\medskip

This finishes the proof of Proposition 3.1.

\bigskip

{\it Proof of Theorem 1.1.}  We apply the Gorenstein-Harada
classification of the finite simple groups of sectional 2-rank at
most four (see [G,p.6] or [Su2,chapter 6, Theorem 8.12]). By
Proposition 3.1,  $G$ has sectional 2-rank at most four and hence
is one of the groups in the Gorenstein-Harada list; the groups
are the following:

$$PSL(n,q), \;  PSU(n,q) \; (n \le 5, \; q \; {\rm odd}),$$
$$G_2(q),\; ^3D_4(q),\;  PSp(4,q) \;  (q \; {\rm odd}), \;
^2G_2(3^{2m+1}),$$
$$PSL(2,8),\; PSL(2,16),\; PSL(3,4),\; PSU(3,4),\; Sz(8),$$
$$\A_l \; (7 \le l \le 11), \; M_i \;(i \le 23),\; J_i \;(i \le 3),
\;Mc,\; Ly$$

\bigskip

The alternating group $\A_7$ has the linear
fractional group  ${\rm PSL}(2,7)$ (which has a subgroup $\S_4$
of index seven) as a subgroup, so by Proposition 2.1 it does not act
on a homology 4-sphere. This excludes $\A_n$ for $n
\ge 7$. Also, the sporadic groups $M_i$, $J_i$, $Mc$ and $Ly$ in
the Gorenstein-Harada list have metacyclic subgroups $H$ excluded by
Lemma 2.2 or linear fractional subgroups
${\rm PSL}(2,q)$ excluded by Proposition 2.1, so they do not act on a
homology 4-sphere (see [A] for information about the maximal
subgroups of the sporadic groups).

\medskip

All other groups in the list
except ${\rm PSL}(2,5)$ and ${\rm PSL}(2,9)$ can be excluded along
similar lines. For example, concerning the linear groups $L_m(q) =
{\rm PSL}(m,q)$,
$q=p^n$, we note that $L_2(p)$ is a subgroup of
$L_2(q)$; also, for $m>r$, the linear group
${\rm SL}(r,q)$ is a subgroup of the linear fractional group
$L_m(q) = {\rm PSL}(m,q)$ (see also [Su2,chapter 6.5] where the
centralizers of involutions in the classical groups are
determined).  Applying Proposition 2.1 and Lemma 2.4, it suffices
then to exclude the groups
$L_3(2)$, $L_3(3)$ and $L_3(5)$. But
$L_3(2)$ is isomorphic to $L_2(7)$, the group
$L_3(3)$ has a metacyclic subgroup
$H(13:3)$ excluded by Lemma 2.2 ([A] or [Su2,p.530]), and
$L_3(5)$ a metacyclic subgroup $H(31:3)$. Thus among the linear
fractional groups there remain only
$L_2(5) \cong \A_5$ and $L_2(9) \cong \A_6$.

\medskip

We will not repeat the arguments for the other groups: the most
interesting remaining case is the Suzuki group
$Sz(8)$ which has one conjugacy class of involutions and a subgroup
$(\Z_2)^3$; applying the Borel formula to this subgroup (similar
as for $PSL(2,8)$ in the proof of Proposition 2.1) shows that
$Sz(8)$ does not act on a homology 4-sphere.

\medskip

This finishes the proof
of Theorem 1.1  (we remark that in fact all finite simple groups
exept $\A_5$ and $\A_6$ can be excluded in the same way).

\bigskip

{\bf 4. Proof of Theorem 1.2}

\medskip

In the proof of Theorem 1.2 we need some extra  informations about the  elementary abelian groups acting on homology $4$-spheres;  we summarize  them in the following lemma. We  recall that the fixed point set of a
group of prime order  $\Bbb{Z}_p$, that acts
orientation-preservingly on a homology $4$-sphere, is a homology
sphere of dimension zero or two ([Bo, chapter IV.4]). To obtain Lemma 4.1. we use again  the Borel Formula but in a more technical way.
\bigskip

{\bf 4.1 Lemma}  {\sl For a prime $p$, let $A$ be an elementary
abelian
$p$-group acting orientation-preservingly on a homology 4-sphere.

a) If  $A$ has rank two and $p$ is odd, then $A$ contains
exactly two cyclic subgroups with 0-dimensional fixed point set.

b) If $A$ has rank two and $p=2$, then $A$ contains at least one
involution with 2-dimensional fixed point set; in the case of
three involutions with 2-dimensional fixed point set the group
$A$ has a fixed-point set of dimension one.

c) If $A$ has rank three and $p=2$, then $A$ can contain either
one or three involutions with 0-dimensional fixed point set.

d) If $A$ has rank four and $p=2$ , then $A$ contains exactly
five involutions with 0-dimensional fixed-point set.}
\bigskip

{\it Proof.} For an  abelian $p$-group $A$  acting on a
homology 4-sphere the Borel Formula  appears as follow:

                $$4-r = \Sigma_H (n(H)-r).$$

The sum is taken over the subgroups
$H$ of index $p$ in  $A$, $n(H)$ denotes the dimension of the
 fixed point set of $H$ and $r$ the dimension of the  fixed point
set of $A$ (equal to -1 if the fixed point set is empty).

If $A$ has rank two  we have $p+1$ cyclic subgroups of index
$p$; if $H$ is cyclic then either $n(H)=0$ or $n(H)=2$.
Considering all the possibilities we obtain the situations
described in points a) and b).  In particular for the 2-groups
we have:  if $r=1$ then $A$ contains three involutions with
2-dimensional fixed point set, if $r=0$  then $A$ contains two
involutions with 2-dimensional fixed point set,    if $r=-1$
then $A$ contains one involution with 2-dimensional fixed point
set (the case $r=2$ can not occur). We remark that the dimension
of the global fixed point set of the group determines  the
fixed point sets of the involutions contained in the group.

So we apply the formula to $A=(\Bbb{Z}_2)^3$; in this case we
have seven subgroups of index two that are elementary abelian
2-groups of rank two. A priori $r$ can be equal to $1,$ $0$ or
$-1$ and the possibilities for
$n(H)$ are given by the previous case.  We explain the situation
for  $r=0$; in this case we obtain from the formula that we can
have either four subgroups with 1-dimensional fixed point set
and three  subgroups with 0-dimensional fixed point set or five
subgroups with 1-dimensional fixed point set, one   subgroup with
0-dimensional fixed point set  and  one   subgroup with empty
fixed point set (it is possible to see also  that the second
case   can not occur but it is not necessary for the proof).
The fixed point set of a subgroup of index two determines the
fixed point sets of the involutions contained in the subgroup;
since any involution is contained in three different subgroups
of index two,  we can compute that in both cases  we have six
involutions with 2-dimensional fixed point set and one
involution with 0-dimensional fixed point set.  Analogously we
can analyze the case $r=-1$; the case $r=1$ can not occur.

For the case of elementary abelian group of rank four we can
repeat a similar computation referring to the results for
elementary abelian group of rank three; this concludes the
proof.

\bigskip

Recall that a finite $Q$ group is {\it quasisimple} if it is
perfect (the abelianized group is trivial) and the factor group
of $Q$ by its center is a non-abelian simple group.  A finite
group $E$ is  {\it semisimple} if it is perfect and the factor
group of $E$ by its center is a direct product of simple
non-abelian groups (see [Su2, chapter 6.6]). A semisimple group
is a central  product of quasisimple groups that are uniquely
determined. Any finite group $G$ contains a  unique maximal
semisimple normal group
$E(G)$ (the subgroup $E(G)$ may be trivial); the subgroup $E(G)$
is characteristic in $G$ and the quasisimple factors of $E(G)$
are called the {\it components} of $G$.     To prove the Theorem
1.2  we consider first the case of groups with trivial maximal
normal semisimple subgroup.

\bigskip

{\bf 4.2  Lemma}  {\sl Let $G$ be a group acting
orientation-preservingly on a homology $m$-sphere such that $G$
has trivial maximal normal semisimple subgroup $E(G)$.

a) If $m=3$ then $G$ is solvable.

b) If $m=4$ either  $G$ is solvable or  $G$ contains a normal
subgroup isomorphic to $(\Bbb{Z}_2)^4$ with factor group
isomorphic  to
$\Bbb{A}_5$ or  $\Bbb{S}_5$.}

\bigskip

{\it Proof.} a) This is  [MeZ, Proposition 8].

b)  We consider first the  case of $G$ containing  a normal
non-trivial cyclic subgroup $H$ and we prove that in this case
if $E(G)$ is trivial then $G$ is solvable. We can suppose that
$H$ has prime order $p$ so each element of $H$ has the same
fixed point-set; since $G$ normalizes $H$ then $G$  fixes
setwise the fixed point set of $H$.

If the fixed point set of $H$ is $S^0$ there exists a subgroup
$G_0$  of index at most two in $G$ such that $G_0$ fixes both
points of $S^0$; the subgroup $G_0$ acts faithfully on a
3-sphere, that is the boundary of a regular invariant
neighbourhood of one of the two fixed points. Since
$E(G_0)$ is trivial  we obtain by a) that $G^0$, and
consequently $G$, are solvable.

If the fixed point set is $S^2$ we consider in $G$  the normal
subgroup
$K$ of elements fixing pointwise $S^2$; the subgroup $K$
contains $H$ and since $K$ acts locally as a rotation around
$S^2$ then $K$ is cyclic. The factor group $G/K$ acts faithfully
on $S^2$. If  $G/K$ is solvable, we get the thesis; otherwise we
can suppose that
$G/K$ is isomorphic either to  $\Bbb{A}_5$ or to $\Bbb{Z}_2\times
\Bbb{A}_5$  because these  are the only non-solvable finite
groups acting on $S^2$ (the action on $S^2$ is not necessarly
orientation-preserving).       In both cases  the action of
$\Bbb{A}_5$  by conjugation  on $K$ is trivial because $K$ is
cyclic and its  automorphism group is abelian; then
$G$ contains with index at most two $G_0$ a subgroup that is a
central extension of $\Bbb{A}_5$. The derived group $G'_0$ is a
quasisimple normal subgroup of $G_0$  (see [Su1, Theorem 9.18,
pag.257]) and this fact implies that  $E(G)$ is not trivial in
contraddiction with our hypothesis.

The proof of this particular case is now complete  and in the
following we can use this fact.

\smallskip

{\it Fact: if a subgroup $N$ of $G$ contains  a non-trivial
cyclic normal subgroup then either  $N$ is solvable or  $E(N)$
is not trivial.}

\smallskip

We consider now the general case.  We denote by $F$  the Fitting
subgroup of $G$ (the maximal nilpotent normal subgroup of $G$).
Since  $E(G)$ is trivial,  the Fitting subgroup
$F$ coincides with   the generalized Fitting subgroup that  is
the product of the Fitting subgroup with the maximal semisimple
normal subgroup. The generalized Fitting subgroup $F$ contains
its centralizer in $G$  and
$F$ is not trivial ([Su2, Theorem 6.11, pag.452]).

Since $F$ is nilpotent it is the direct product of its Sylow
$p$-subgroups. In particular any Sylow subgroup of $F$ is normal
in $G$; since $F$ is not trivial we have  $P$ a non-trivial
$p$-subgroup normal in $G$. We consider $Z$ the maximal
elementary abelian $p$-subgroup contained in the center of $P$;
this subgroup is not trivial and it is normal in $G$.

Suppose first that we can chose $p$  odd (the order of $F$ is
not a power of two); then, by Lemma 2.3, $Z$ has
rank one or two.  If $Z$ is cyclic, by the first part of the
proof,   $G$ is solvable. If $Z$ has rank two by Lemma 4.1 it
contains exactly two cyclic subgroups
$H$ and $H'$ with 0-dimensional fixed point set; $G$ acts by
conjugation on the set $\{H,\,H'\}$ and $N_G{H}$ has index at
most two in $G$. Since $E(G)$ is trivial then the maximal normal
semisimple subgroup of
$N_G(H)$ is trivial; we obtain that $N_G(H)$, and consequently
$G$, are solvable.  This concludes this case.

Suppose now that the order of $F$ is a power of two; in this
case $F=P$ is a 2-group and $Z$ is an elementary abelian 2-group
of  rank at most four (by Lemma 2.3).

If $Z$ has rank one by the first part of the proof $G$ is
solvable.

If $Z$ has rank two  we consider $C_G(Z)$ the centralizer of $Z$
in $G$ that is normal bacause $Z$ is normal;   $C_G(Z)$ contains
a non-trivial normal cyclic subgroup and it is solvable. The
factor $G/C_G(Z)$ is isomorphic to  a subgroup of $GL(2,2)$, the
automorphism group of an elementary abelian 2-group of rank two;
since  $GL(2,2)$ is a solvable group we can conclude that $G$ is
solvable.

Suppose that  $Z$ has rank three.  In this case the factor group
$G/C_G(Z)$ is  isomorphic to  a subgroup of $GL(2,3)$, the
automorphism group of an elementary abelian 2-group of rank
trhee; $GL(2,3)$ has order
$2^3\cdot 3\cdot 7$ and any element of order seven permutes
cyclically all the involutions of $(\Bbb{Z}_2)^3$. The group
$G/C_G(Z)$ can not contain element of order $7$ otherwise all
involutions in $Z$  are conjugated and this is impossible by
Lemma 4.1; so the group $G/C_G(Z)$ has order $2^\alpha 3^\beta$
and it is solvable. This fact implies that  $G$ is solvable.

It remains the case $Z$ of rank four; by Lemma 4.1 the group
$Z$ contains at least one involution  with  fixed point set
$S^0$. The group $P$ fixes setwise $S^0$ and we have a subgroup
$P_0$ of index at most two that fixes both the points in $S^0$.
The center of $P_0$ contains  an elementary abelian group of
rank at least three (the group $Z\cap P_0$) and $P_0$ acts
faithfully on a 3-sphere that is the boundary of a regular
invariant neighbourhood of one of the two fixed points;  the
only 2-group acting on the 3-sphere  with this property is
$(\Bbb{Z}_2)^3$ (see [MeZ, Proposition 2 and Proposition 3]).
In this case the generalized Fitting subgroup $F=P$ is an
elementary abelian group of rank four, we recall that $F$
contains its centralizer in
$G$. By Lemma 4.1 in $F$ are contained exactly five involutions
with 0-dimensional fixed point set; they generate the group $F$
because by Lemma 4.1 the  subgroups of index two  contain at
most three of such involutions. We can conclude that $G/F$ acts
faithfully on the set of the  five involutions of $F$ with
0-dimensional fixed point set and $G/F$ is isomorphic to a
subgroup of $\Bbb{S}_5$.  In particular if $G$ is not solvable
we obtain that $G/F$ is isomorphic either to $\Bbb{A}_5$ or to
$\Bbb{S}_5$ (the only non-solvable subgroups of  $\Bbb{S}_5$). This
fact concludes the proof.

\bigskip

Now we consider the case of semisimple group.

\bigskip

{\bf 4.3  Lemma}  {\sl Let $G$ be a finite semisimple group
acting orientation-preservingly on a homology $4$-sphere, then
$G$ is isomorphic to one of the following group:
$$\Bbb{A}_5,\,   \Bbb{A}_6,\,
\Bbb{A}_5^*,\,\Bbb{A}_5^*\times_{\Bbb{Z}_2}\Bbb{A}_5^*.$$ }

{\it
Proof.} By Theorem 1.1 and Lemma 2.4,   if $G$ is quasisimple,
then $G$  is isomorphic to $\Bbb{A}_5$, $ \Bbb{A}_6$ or to
$\Bbb{A}_5^*\cong$SL(2,5) that is the unique perfect central
extension of $\Bbb{A}_5$.

We consider now the case of   $G$ with two quasisimple
components; since in our list of quasisimple groups
$\Bbb{A}_5^*$  is the unique group with non-trivial center, then
either $G\cong
\Bbb{A}_5^*\times_{\Bbb{Z}_2}\Bbb{A}_5^*$ or $G$ is the direct
product of two quasisimple subgroups. We have to exclude this
second possibility  thus we consider $G\cong Q\times Q'$ with
$Q$ and $Q'$ isomorphic to  $\Bbb{A}_5$, $ \Bbb{A}_6$ or
$\Bbb{A}_5^*$.

Suppose first that one of the two components, say $Q$, is
isomorphic to
$\Bbb{A}_6$. Let $f$ be a non trivial  element of $Q'$, since
$Q$ commutes elementwise with $f$, the subgroup $Q$ fixes
setwise the fixed point set of $f$.

If the fixed point set of $f$ is a 0-sphere  there exists a
subgroup of
$Q$ of index at most two that acts  faithfully on a 3-sphere
that is the boundary of a regular invariant  neighbourhood of
one of the two fixed points. The group  $\Bbb{A}_6$ has no
subgroup of index two and  $\Bbb{A}_6$ can not  act faithfully
on a homology 3-sphere (see [MeZ, Theorem 2]).

If the fixed point set of $f$ is a 2-sphere there exists a
factor group of
$Q$ by a cyclic subgroup that acts faithfully on a 2-sphere.
Since $\Bbb{A}_6$ is simple and $\Bbb{A}_6$ can not act on $S^2$
also this case can not occur.

Suppose now that one component is isomorphic to $\Bbb{A}_5^*$,
in this case the center of $G$ is not trivial because it
contains the involution that is in  the center of
$\Bbb{A}_5^*$.  The same argument used for $\Bbb{A}_6$ applies
to $G$ and we can exclude
$\Bbb{A}_5^*$ as component.

The only case that remains is $G\cong \Bbb{A}_5\times
\Bbb{A}_5$. We denote by $A$  the Sylow 2-subgroup of the first
component. The subgroup
$A$ is elementary abelian of rank two and the three involutions
$\{t_1,t_2,t_3\}$ in $A$ are all conjugated; by Lemma 4.1  the
only possibilities is that the fixed point set of each
involution in $A$ is a 2-sphere. We consider  the  group
generated by $t_1$ and by the second component; this group is
isomorphic to $\Bbb{Z}_2\times \Bbb{A}_5$ and  acts faithfully
on the  2-sphere that is the fixed point set of $t_2$. Since
$t_1$ commutes with $A_5$,  the action of $t_1$ on the fixed
point set of $t_2$ is free in contradiction with  Lemma 4.1.

So we have that, if $G$ has two components, $G$ is isomorphic to
$
\Bbb{A}_5^*\times_{\Bbb{Z}_2}\Bbb{A}_5^*$.

Finally we exclude the possibility to have more then two
components. By the previous cases we can argue  that all the
components are isomorphic to  $ \Bbb{A}_5^*$, but in this case
the center of $G$ is not trivial and referring again either to
three dimensional case  ([MeZ, Theorem 2)] or to  two
dimensional case we can exlcude these groups.

\bigskip

{\it Proof of Theorem 1.2} If the maximal semisimple
normal subgroup $E$ of $G$ is trivial we apply Lemma 4.2 and we
obtain that either $G$ is solvable or we are in case a).

Suppose then that $E$ is not trivial, then $E$ is one of the
groups presented in Lemma 4.1. Since $E$ is normal, its
centralizer $C=C_G(E)$ in
$G$ is also normal. We denote by $\tilde E$ the subgroup
generated by  $C$ and $E$; by definition $\tilde E$ is a  normal
subgroup of $G$  that  is a central product of $C$ and $E$. The
factor group $G/\tilde E$ is  a subgroup of the outer
automorphism group Aut$(E)/$Inn$(E)$ of $E$.

The maximal semisimple normal subgroup of $C$ is trivial
otherwise $E$ is not maximal in $G$. So by Lemma 4.2 either $C$
is solvable or $C$ contains an elementary normal  2-subgroup of
rank four and $C$ has trivial center; the second case can not
occur otherwise $\tilde E$ would contain an elementary abelian
group of rank five in contradiction with Lemma 2.3.

\medskip

We consider  first the case $E\cong \Bbb{A}_6$; the same arguments used
in Lemma 4.3 to exclude $ \Bbb{A}_6\times Q$ apply and we
obtain that $C$ is trivial. The outer automorphism group of
$\Bbb{A}_6$ is isomorphic to $\Z_2 \times \Z_2$ (see [Su1,p.300]), so
there are three 2-fold extensions of $\A_6 \cong {\rm PSL}(2,9)$
corresponding to the three subgroups $\Z_2$ of $\Z_2 \times \Z_2$. One
of these extensions is $\S_6$ which acts orientation-preservingly and
orthogonally on the 4-sphere. Another 2-fold extension is ${\rm
PGL}(2,9)$. The Sylow 3-subgroup of ${\rm PGL}(2,9)$ is isomorphic to
$\Z_3 \times \Z_3$.  In ${\rm PGL}(2,9)$ all elements of order
three are conjugate so  their fixed point sets have the same dimension but this is impossible by Lemma 4.1.
 Hence ${\rm PGL}(2,9)$ does not act on a homology
4-sphere. A similar argument works for the third 2-fold extension of
${\rm PSL}(2,9)$ (the Matthieu group $M_{10}$, see [A]) in which also
all elements of order three are conjugate. Now also the
$\Z_2 \times \Z_2$-extension of ${\rm PSL}(2,9)$ does not act, and we
are left with the two groups $\A_6$ and $\S_6$ in case b) of Theorem
2.1.

\medskip

Suppose now that  $E\cong \Bbb{A}_5^*\times_{\Bbb{Z}_2}\Bbb{A}_5^*
$. We denote  by $t$ the central involution in $E$ that is
central in   $G$ because $E$ is characteristic in $G$. The fixed
point set of $t$ is not a 2-sphere otherwise a factor of
$\Bbb{A}_5^*\times_{\Bbb{Z}_2}\Bbb{A}_5^* $  by a cyclic group
would act faithfully on the 2-sphere.  The fixed point set of
$t$ is $S^0$  and we get $G_0$ a subgroup of $G$ of index at
most two acting faithfully and orientation presenvingly  on a
3-sphere, that is the boundary of a regular invariant
neighbourhood of one of the two fixed points. Since
$\Bbb{A}_5^*\times_{\Bbb{Z}_2}\Bbb{A}_5^* $ is maximal beetwen
the groups acting on the  3-sphere (see [MeZ, Theorem 2]) we
obtain that
$G_0\cong \Bbb{A}_5^*\times_{\Bbb{Z}_2}\Bbb{A}_5^*$.

\medskip

Next we consider the case of $E\cong \A_5$. A Sylow 2-subgroup
$S$ of $E$ is an elementary abelian group of rank two and the
three involutions
$\{t_1,t_2,t_3\}$ in $S$ are all conjugated; by Lemma 4.1 we
have that all the three involutions have 2-dimensional fixed
point set. We consider
$S^2$ the fixed point set of $t_1$, the action of $t_2$ on $S^2$
is orientation reversing  and then $C$ acts
orientation-preservingly and faithfully on $S^2$. Suppose $C$
isomorphic to $ \Bbb{A}_5$,  $ \Bbb{S}_4$ or $ \Bbb{A}_4$, since
$t_2$ commutes with $C$,  then $t_2$ has to act freely on $S^2$
in contraddiction with Lemma 4.1. We can conclude that the group
$C$ can be either dihedral  or cyclic. The outer automorphism
group of $\Bbb{A}_5$ has order two and the factor group
$G/\tilde E$ has order at most two.

\medskip

Suppose finally  that $E\cong \Bbb{A}_5^*$; we denote by $t$ the
central involution in $E$,  $t$ is central also in $G$ since $E$
is characteristic. Suppose that the fixed point set of $t$ is
$S^0$, a 0-sphere. In this case we get $G_0$ a subgroup of $G$
with index at most two, that acts faithfully and orientation
preservingly  on a 3-sphere. By [MeZ, Theorem 2] we can conclude
that $G_0$ is isomorphic to $G_0\cong
\Bbb{A}_5^*\times_{\Bbb{Z}_2} C$ where $C$ is a solvable group
acting freely on a 3-sphere (for further details see also [Z2]).

If the fixed point set of $t$ is $S^2$ we denote by $K$ the
cyclic group of elements that fix pointwise $S^2$. The factor
$G/K$  acts faithfully on
$S^2$, then $G/K$ contains with index at most two $\Bbb{A}_5$;
since $K$ is cyclic the action by conjugation of $\A_5$ is
trivial on $K$.    In this case $G$ contains a subgroup of index
at most two isomorphic to
$\Bbb{A}_5^*\times_{\Bbb{Z}_2} \Bbb{Z}_{2n}$. Any cyclic group
admits  a free and orientation preserving action on the 3-sphere
(but in this case the action is not naturally related to the
action on the homology 4-sphere).

\medskip

We have considered all the possible cases for $E$ and the proof
is finished.

\bigskip \bigskip

\centerline {\bf References}

\bigskip

\item {[A]} J.H.Conway, R.T.Curtis, S.P.Norton, R.A.Parker,
R.A.Wilson, {\it Atlas of Finite Groups.} Oxford University
Press 1985

\item {[Bo]} A.Borel, {\it Seminar on Transformation Groups.}
Annals of Mathematics Studies 46, Princeton University Press 1960

\item {[Br]} G.Bredon, {\it Introduction to Compact
Transformation Groups.} Academic Press, New York 1972

\item {[DH]}  R.M.Dotzel, G.C.Hamrick, {\it  $p$-group actions
on homology spheres.}  Invent. math. 62, 437-442  (1981)

\item {[DV]}  P.Du Val, {\it  Homographies, Quaternions and Rotations.}
Oxford Math. Monographs, Oxford University Press 1964

\item {[G]} D.Gorenstein, {\it The Classification of Finite
Simple Groups.} Plenum Press, New York 1983

\item {[MeZ]} M.Mecchia, B.Zimmermann, {\it On finite groups
acting on $\Z_2$-homology 3-spheres.}
 Math. Z. 248, 675-693 (2004)

\item {[Mn]} J.Milnor, {\it Groups which act on $S^n$ without
fixed points.} Amer. J. Math. 79, 623-630  (1957)

\item {[Sj]} D.Sjerve, {\it Homology spheres which are covered by
spheres.}  J. London Math. Soc. 6, 333-336  (1973)

\item {[Su1]} M.Suzuki, {\it Group Theory I.}  Springer-Verlag
1982

\item {[Su2]} M.Suzuki, {\it Group Theory II.}  Springer-Verlag
1982

\item {[Sw]} R.G.Swan, {\it The $p$-period of a finite group.}
Illinois J. Math. 4, 341-346 (1960)

\item {[Z1]} B.Zimmermann, {\it On finite simple groups acting
on homology 3-spheres.}  Topology Appl. 125, 199-202 (2002)

\item {[Z2]} B.Zimmermann, {\it On the classification of finite
groups acting on homology 3-spheres.}
Pacific J. Math. 217, 387-395 (2004)

\item {[Z3]} B.Zimmermann, {\it Cyclic branched coverings and
homology 3-spheres with large group actions.}
Fund. Math. 184, 343-353  (2004)

\bye